\documentclass[11pt]{amsart}
\usepackage{amssymb,amsmath}

\newcommand\Ella[1]{\ell(#1)}
\newcommand\Qlla[2]{(-q)^{{#1}\ell(#2)}}
\newcommand{\II}{\mathbb I}
\newcommand{\K}{\Bbbk} 
\newcommand{\Sym}{\mathfrak S}
\newcommand{\Hom}{\mathrm{Hom}}
\newcommand{\newtilde}[1]{\tilde{#1}} 
\newcommand{\qdet}{\mathrm{det}_q\,}

\theoremstyle{plain} 
\newtheorem{mytheo}{Theorem}
\newtheorem{myprop}[mytheo]{Proposition}

\theoremstyle{definition} 
\newtheorem{mydefi}{Definition}

\theoremstyle{remark} 
\newtheorem*{mynota}{Notation} 

\newtheorem*{myrema}{Remark}
\newtheorem*{myexam}{Example}

\def\TakNote{Compare the discussion of $N(W)$ in loc. cit., especially Proposition 3.3.}
\title{%
	A class of left quantum groups 
	modeled after $\mathrm{SL}_q(n)$}

\author{Aaron Lauve}
\address[Lauve]{%
	Department of Mathematics\\ 
        Texas A\&M University\\
        College Station, TX \ 77843
         }
\email{lauve@math.tamu.edu}
  \urladdr{http://www.math.tamu.edu/\~{}lauve}

\author{Earl J. Taft}
\address[Taft]{%
	Department of Mathematics\\
        Rutgers, The State University of New Jersey\\
        110 Fre\-ling\-huysen Road\\
        Piscatway NJ, 08854-8019, USA}
\email{etaft@math.rutgers.edu}

\begin{document}

\date{13 March 2006}

\begin{abstract}
For each $n\geq 2$, we construct a left quantum group, i.e., a left Hopf algebra $\newtilde{\mathrm{SL}}_q(n)$ generated by comatrix units $X_{ij}$, $1\leq i,j\leq n$, which has a left antipode but no  right antipode. The quantum special linear group $\mathrm{SL}_q(n)$ is a homomorphic image of $\newtilde{\mathrm{SL}}_q(n)$.
\end{abstract}

\maketitle

\section*{Notation}
\noindent We collect a few standard and non-standard notations:

Let $[n]$ denote the set $\{1,2,\ldots,n\}$, and let $[n]^k$ denote the set of $k$-tuples chosen from $[n]$. Given a $k$-tuple $I = (i_1,i_2, \ldots, i_k)$, let $\Ella{I}$ denote its \emph{length}, i.e., the least number of adjacent interchanges necessary to put the elements of $I$ in non-decreasing order. Define $\ell$ on permutations $\pi\in\Sym_k$ via the standard one line notation, i.e., $\ell(\pi) = \ell(\pi(1),\pi(2),\ldots, \pi(k))$; also, denote $\pi(k)$ by $\pi_k$.

For an $n$-tuple $I=(i_1,\ldots,i_n)$, say $I\in\Sym_n$ if the elements $i_1,\ldots, i_n$ are distinct, and $I\not\in\Sym_n$ otherwise. We use the \emph{membership} Kronecker delta function $\delta_{I,\Sym_n}$ to distinguish the two cases, taking value one in the former case and zero in the latter case. 

Throughout,  $\K$ will denote a field containing a distinguished invertible element $q$. For any $\K$-vector space $V$, we let $\II$ denote the identity mapping on $V$; similarly, given a $\K$-algebra $R$, we let $\II_n$ be the identity matrix in $M_n(R)$.

\section{Introduction}\label{sec:intro}

Let $B$ be a bialgebra over $\K$. $B$ has an associative multiplication $m:B\otimes B\rightarrow B$, a unit element $\mu: \K \rightarrow B$, a coassociative comultiplication $\Delta: B\rightarrow B\otimes B$, and a counit $\varepsilon : B\rightarrow \K$ such that $\Delta$ and $\varepsilon$ are $\K$-algebra maps. The $\K$-linear maps $\Hom(B,B)$ form a monoid with respect to the convolution product $f\ast g = m(f\otimes g)\Delta$ and the unit element $\mu\varepsilon$.

A bialgebra $H$ is a Hopf algebra if the identity map $\II\in \Hom(H,H)$ is invertible. That is, if there is an \emph{antipode} $S\in\Hom(H,H)$ satisfying $S\ast \II = \II\ast S = \mu\varepsilon$ as functions on $H$. Explicitly, if $\Delta h = \sum_i h_{i} \otimes h'_{i}$, 
$S$ must satisfy
\[
(\forall h\in H) : \sum_i S(h_i)h'_i = \varepsilon(h)1 = \sum_i h_i S(h'_i).
\]
Such an $S$ is unique, and is an algebra and coalgebra antimorphism of $H$.

A bialgebra $H$ is called a \emph{left Hopf algebra} if there is an $S$ in $\Hom(H,H)$ such that $S\ast \II = \mu\varepsilon$; such an $S$ is called a \emph{left antipode} of $H$. If $S$ is not also a right antipode, i.e. $\II\ast S \neq \mu\varepsilon$, then $H$ will have an infinite number of left antipodes \cite{Jac:3}.

The first examples of left Hopf algebras which are not Hopf algebras were constructed in \cite{GreNicTaf:1}. These were free left Hopf algebras on the coalgebras of $n\times n$ comatrices for $n\geq 2$. The specific left antipode constructed was both an algebra and a coalgebra antimorphism.

A variation of the examples of \cite{GreNicTaf:1} was given in \cite{NicTaf:1}. Both constructions begin by defining $S$ on a set of algebra generators. The key change is that one extends $S$ to the whole algebra not by assuming an algebra antimorphism property, but by defining it directly on a basis for the algebra. One obtains a basis of irreducible words in the generators, i.e., one orders the words (monomials) in the generators and uses the Diamond Lemma \cite{Ber:1} to obtain a basis of irreducible words. Then one defines $S$ on an irreducible word by reversing its order on the generators. Roughly speaking, $S$ is ``locally'' an algebra antimorphism, but not ``globally.'' For the examples in \cite{NicTaf:1}, it turns out that no left antipode is an algebra antimorphism (nor a coalgebra antimorphism).

With the current interest in quantum groups, it would be of interest to construct a one-sided quantum group, say a left quantum group which is not a right quantum group; e.g., a left Hopf algebra generated by comatrix units $X_{ij}$, $1\leq i,j\leq n$ for $n\geq 2$, which satisfies some relations connected to those of a known (two-sided) quantum group, but which is not a (two-sided) Hopf algebra. We hope that this would be of interest in quantum physics. We now indicate one possible connection, namely the \emph{boson-fermion correspondence.}

The first attempt to construct a one-sided quantum group was in \cite{RodTaf:1}. There one started with 3 of the 6 relations for a $2 \times 2$ quantum matrix $X$, namely
\[
\tag{$\dagger$} \left.
\begin{array}{ll}
X_{21}X_{11} &= qX_{11}X_{21} \\
X_{22}X_{12} &= qX_{12}X_{22} \\
X_{22}X_{11} &=X_{11}X_{22}- q^{-1}X_{21}X_{12}+qX_{12}X_{21}
\end{array}\right\}
\]
(see relations (1)-(3) of \cite{RodTaf:1}). In \cite{RodTaf:1}, requiring the left antipode $S$ to be an algebra antimorphism necessitated additional relations, and the result was a two-sided
antipode. Thus while the result was a new quantum group, it was not an example of a one-sided quantum group. 

The relations ($\dagger$) have been generalized in \cite{GarLeZei:1} to $r \times r$ matrices $A=(a_{ij})$ by making every $2 \times 2$ submatrix of $A$ satisfy ($\dagger$). The resulting algebra is called the right quantum algebra and $A$ is called a right quantum matrix.  This one-sided quantum setting is sufficient to prove the main Theorem 1 of \cite{GarLeZei:1}, a \emph{quantum} boson-fermion correspondence.

Write $\mathrm{Bos}(A) = \sum_{(m_1, \ldots, m_r)\in \mathbb{N}^r} G(m_1, \ldots, m_r)$, where $G(m_1, \ldots, m_r)$ is the coefficient of $x_1^{m_1}\cdots x_r^{m_r}$ in the product $X_1^{m_1} \cdots X_r^{m_r}$ (taking $X_i := \sum_{j=1}^{r} a_{ij}x_j$ and the $x_i$'s to be commuting variables which commute with the $a_{jk}$'s). Write $\mathrm{Ferm}(A) = \sum_{J\subseteq[r]} (-1)^{|J|}\qdet A_J$, where $A_J$ is the $J\times J$ submatrix of $A$ and $\qdet A_J$ is the usual quantum determinant of quantum group theory. Theorem 1 of \cite{GarLeZei:1} states that $\mathrm{Bos}(A) = 1/\mathrm{Ferm}(A)$. Another proof appears in \cite{FoaHan:2}. 

When $q=1$ and $A$ has commuting entries, $\sum_{m_1+\cdots + m_r=n} G(m_1, \ldots, m_r)$ is the trace of $S^n(A)$, the $n$-th symmetric power of $A$, and $\det (I-tA) = \sum_{n\in\mathbb N} (-1)^n \mathrm{tr}\, \Lambda^n(A)t^n$, where $\Lambda^n$ is the $n$-th exterior power of $A$. The boson-fermion correspondence states that $\sum_n \mathrm{tr}\, S^n(A)t^n$ and $\sum_n (-1)^n \mathrm{tr}\, \Lambda^n(A)$ are inverses of each other. Thus the main result of \cite{GarLeZei:1} is a quantum generalization of the boson-fermion correspondence.

In \cite{RodTaf:2}, the third relation of ($\dagger$) was split into two relations, each setting a version of the quantum determinant equal to one. A basis of irreducible words was obtained, and the ideas of \cite{NicTaf:1} were used to obtain a left (but
not right) quantum group. This construction is reviewed in Section 2.
Here, we generalize the construction in \cite{RodTaf:2} to produce a left quantum group $\newtilde{\mathrm{SL}}_q(n)$ for all $n\geq 2$. The result for $n>2$ is an algebra which possesses $\mathrm{SL}_q(n)$ as a quotient, but does not satisfy any of the standard quadratic relations of $\mathrm{SL}_q(n)$. Rather, there are $n^n$ relations defining $\newtilde{\mathrm{SL}}_q(n)$: $n!$ of them are of the form saying a certain version of the quantum determinant of $X$ is equal to one; the other $n^n-n!$ relations are homogeneous of degree $n$. As for $n=2$, we define the left antipode on generators using a quantum adjoint matrix, then extend $S$ to $\newtilde{\mathrm{SL}}_q(n)$ using the idea of \cite{NicTaf:1}.

\section{The Case $n=2$}\label{sec:nequals2}

We review some elements of the construction in \cite{RodTaf:2}. Consider the bialgebra $\K\langle X\rangle = \K\langle X_{11}, X_{12}, X_{21}, X_{22}\rangle$ with $\Delta X_{ij} = \sum_{k=1}^2 X_{ik}\otimes X_{kj}$ and $\varepsilon(X_{ij}) = \delta_{ij}$ for $1\leq i,j\leq 2$. Ordering the $X_{ij}$ lexicographically on $(i,j)$, we impose four relations, viewed as reduction formulas:
\begin{align}
X_{21} X_{11} &= q X_{11} X_{21} \label{eq:SL2-1} \\
X_{22} X_{12} &= q X_{12} X_{22} \label{eq:SL2-2} \\
X_{22} X_{11} &= q X_{12} X_{21} + 1 \label{eq:SL2-3} \\
X_{21} X_{12} &= q X_{11} X_{22} - q1 \label{eq:SL2-4} 
\end{align}
Applying the Diamond Lemma, one gets a basis for $\newtilde{\mathrm{SL}}_q(2)$ comprising irreducible words which behave nicely under comultiplication---our ``{property ($\star$)}'' in Section \ref{sec:mainresults}. From here it is straightforward to define a left but not a right antipode $S$. 
As happens with the examples in \cite{NicTaf:1}, $\newtilde{\mathrm{SL}}_q(2)$ has the property that no left antipode is an algebra antimorphism. Nor is the given $S$ a coalgebra antimorphism (e.g., the condition fails on $X_{11}^{\phantom{1}\,2}$).

Relations (\ref{eq:SL2-1})--(\ref{eq:SL2-4}) hold in $\mathrm{SL}_q(2)$, and account for roughly half the number of relations necessary to define it. One might expect a similar phenomenon for $n>2$, i.e. that roughly half the relations of $\mathrm{SL}_q(n)$ are needed to build a one-sided Hopf algebra. This approach does not seem to work.

The key for generalizing to $n>2$ lies in a certain $\newtilde{\mathrm{SL}}_q(2)$-comodule. Consider the \emph{quantum exterior plane} $\Lambda_q(2)$, the $\K$-algebra with generators $\xi_1 \,,\, \xi_2$ and relations $\xi_i\xi_i=0$ and $\xi_2\xi_1 = -q^{-1}\xi_1\xi_2$.  In \cite{RodTaf:2} it is shown that $\Lambda_q(2)$ is a right $\newtilde{\mathrm{SL}}_q(2)$-comodule algebra under the mapping $\rho(\xi_i) = \sum_{j=1, 2} \xi_j \otimes X_{ji}$. In fact, this \emph{explains} the defining relations for $\newtilde{\mathrm{SL}}_q(2)$. To illustrate, consider the action of $\rho$ on $\xi_1\xi_1$. On the one hand it must be zero (since $\xi_1\xi_1$ is zero). On the other hand, it is the product $\rho(\xi_1)\rho(\xi_1) = \sum_{i,j=1}^2 \xi_i\xi_j\otimes X_{i1}X_{j1} = \xi_1\xi_2\otimes \left(X_{11}X_{21} + (-q^{-1})X_{21}X_{11} \right)$, from which we deduce 
\begin{align}\tag{$1'$} 
\label{eq:altSL2-1}X_{11}X_{21} - q^{-1} X_{21}X_{11} = 0.
\end{align} 
Applying $\rho$ to the other monomials of degree two yields
\begin{align}
\label{eq:altSL2-2}\tag{$2'$} X_{12}X_{22} - q^{-1} X_{22}X_{12} &= 0 \\
\label{eq:altSL2-3}\tag{$3'$} X_{12}X_{21} - q^{-1} X_{22}X_{11} &= - q^{-1} \\
\label{eq:altSL2-4}\tag{$4'$} X_{11}X_{22} - q^{-1} X_{21}X_{12} &= 1
\end{align}
(the last two coming from demanding that the group-like element $\qdet X$ take the value one).
Equations ($1'$)--($4'$) are readily generalized to $n>2$.

\section{The General Case $n > 2$}\label{sec:general-n}

Let $\K\langle X\rangle$ be the free matrix bialgebra in $n^2$ indeterminants $X=(X_{ij})$, i.e. the free algebra $\K\langle X\rangle$ with $\Delta(X_{ik}) = \sum_j X_{ij}\otimes X_{jk}$ and $\varepsilon(X_{ik}) = \delta_{ik}$. For each $I=(i_1,\ldots,i_n)\in[n]^n$, let $D_I$ denote the sum
\begin{equation*}\label{eq:def-DI}
D_I = \sum_{\pi\in\Sym_n} \Qlla{-}{\pi} X_{\pi_1 i_1}X_{\pi_2 i_2}\cdots X_{\pi_n i_n}.
\end{equation*}
In the case $I=(1,2,\ldots,n)$ we drop the subscript and just write $D$ for $D_{(1,2,\ldots,n)}$. We call $D$ the \emph{quantum determinant.} These are the left-hand sides of equations ($1'$)--($4'$), generalized to $n>2$. 

\begin{mydefi} Let $\newtilde{\Lambda}_q(n)$ be the $\K$-algebra with generators 
$z_i$ for $1\leq i\leq n$ and relations $\mathcal E_{I}$ for $I\in[n]^n$ defined as follows:
\begin{align*}
\label{eq:A0n-relations}\tag{$\mathcal E_I$}%
z_{i_1}z_{i_2}\cdots z_{i_n} &= \left\{\begin{array}{ll}
\Qlla{-}{I} z_1z_2\cdots z_n & \hbox{if } I=(i_1,\ldots,i_n)\in\Sym_n, \\
0 & \hbox{otherwise.} \\
\end{array}\right.
\end{align*}
\end{mydefi}

This may be viewed as a subspace of the usual quantum exterior space. The next definition comprises the replacements for the usual quantum bialgebra $\mathrm{M}_q(n)$ and Hopf algebra $\mathrm{SL}_q(n)$.

\begin{mydefi} The algebra $\newtilde{\mathrm{M}}_q(n)$ is the quotient of $\K\langle X\rangle$ by the two-sided ideal $E$ generated by elements $E_I$ for $I=(i_1,\ldots,i_n)\in[n]^n$ defined by
\begin{equation}\label{eq:def-EI}
E_I = D_I - \delta_{I,\Sym_n} \Qlla{-}{I} D. 
\end{equation}
The algebra $\newtilde{\mathrm{SL}}_q(n)$ is the quotient of $\newtilde{\mathrm{M}}_q(n)$ by the principal ideal $(D-1)$.
\end{mydefi}

\begin{myrema}
Note that when $I=(1,2,\ldots,n)$, $E_I$ reads $D-D$ so there are $n^n-1$ relations defining $\newtilde{\mathrm{M}}_q(n)$ and $n^n$ relations defining $\newtilde{\mathrm{SL}}_q(n)$. The images of the $X_{ij}$ in the quotients form a linearly independent set of generators, so we abuse notation and write $X_{ij}$ for the image of $X_{ij}$ in the quotients. Similarly, we use the notation $D_I$ for the images of the $D_I$ defined above.
\end{myrema}

$\newtilde{\mathrm{M}}_q(n)$ has its relations defined in such a way as to make $\newtilde{\Lambda}_q(n)$ a right comodule algebra, i.e., the right comodule structure map is an algebra homomorphism compatible with the relations $\mathcal E_I$. Arguing abstractly as in \cite{Tak:3}, one can show that $E$ is the unique minimal (bi)ideal making $\newtilde{\Lambda}_q(n)$ into a $\K\langle X\rangle / E$-comodule algebra.\footnote{\TakNote} However, we consider it instructive, since it reveals the combinatorics in our situation, to show explicit calculations here (cf. Proposition \ref{thm:Mq-comodule}). 

\begin{myexam}\label{exm:SLq3} In the case $n=3$, sending the $E_I$ to zero (and $D$ to $1$) gives 27 relations for $\newtilde{\mathrm{SL}}_q(3)$ of the form $D_{(ijk)}=\Qlla{-}{ijk} \cdot \delta_{(ijk),\Sym_3}$. Among these are $6=3!$ relations saying some version of the quantum determinant equals one. For example, the two relations 
\begin{align}
\nonumber {X_{12}X_{21}X_{32}} - {q^{-1}}{X_{12}X_{31}X_{22}} - {q^{-1}}{X_{22}X_{11}X_{32}} + \vphantom{X} \\
\label{eq:altz2z1z2}\tag{${212}'$}
{q^{-2}}{X_{22}X_{31}X_{12}} + {q^{-2}}{X_{32}X_{11}X_{22}}  - {q^{-3}}{X_{32}X_{21}X_{12}} & = 0 \hphantom{-q^{-3}}
\end{align}
and
\begin{align}
\nonumber {X_{11}X_{23}X_{32}} - {q^{-1}}{X_{11}X_{33}X_{22}} - {q^{-1}}{X_{21}X_{13}X_{32}} + \vphantom{X} \\
\label{eq:altz1z3z2}\tag{${132}'$}
 {q^{-2}}{X_{21}X_{33}X_{12}} + {q^{-2}}{X_{31}X_{13}X_{22}}  - {q^{-3}}{X_{31}X_{23}X_{12}} & = -q^{-1} \hphantom{0}
\end{align}
may be compared with (\ref{eq:altSL2-2}) and (\ref{eq:altSL2-3}) respectively to see the general picture. 
\end{myexam}

We view the $n^n$ relations defining $\newtilde{\mathrm{SL}}_q(n)$ as \emph{reduction formulas} in the proof of Theorem \ref{thm:SLqn-properties} and apply Bergman's Diamond Lemma. For this, we need an ordering on all words in the generators compatible with multiplication. First, we order the generators in lexicographic order according to their subscripts; then for two words $w_1,w_2$ in the generators, say $w_1<w_2$ if it is lower in the $\mathrm{length}+\mathrm{lexicographic}$ order. 

We return to Example \ref{exm:SLq3} for illustration. Put $X_{11} < X_{12} < X_{13} < X_{21} < \cdots < X_{33}$, and order the words in these generators as just described. Written as reduction formulas, equations (\ref{eq:altz2z1z2}) and (\ref{eq:altz1z3z2}) become
\begin{align*}
{X_{32}X_{21}X_{12}} = & {q^3}{X_{12}X_{21}X_{32}} - {q^2}{X_{12}X_{31}X_{22}} 
- {q^2}{X_{22}X_{11}X_{32}} + \hphantom{XXXX} \\
\tag{${212}$}
& {q}{X_{22}X_{31}X_{12}} + {q}{X_{32}X_{11}X_{22}}
\end{align*}
and
\begin{align*}
{X_{31}X_{23}X_{12}} = & {q^3}{X_{11}X_{23}X_{32}} - {q^2}{X_{11}X_{33}X_{22}} 
- {q^2}{X_{21}X_{13}X_{32}} + \hphantom{XXXX} \\
\tag{${132}$}
& {q}{X_{21}X_{33}X_{12}} + {q}{X_{31}X_{13}X_{22}} + {q^2}
\end{align*}
respectively (compare with (\ref{eq:SL2-2}) and (\ref{eq:SL2-3})).

\section{Main Results}\label{sec:mainresults}

\begin{myprop}\label{thm:Mq-comodule} For $\newtilde{\Lambda}_q(n)$ and $\newtilde{\mathrm{M}}_q(n)$, we have:

\noindent \textit{1.} $\newtilde{\mathrm{M}}_q(n)$ is a bialgebra with $\Delta$ and $\varepsilon$ given by $\Delta X_{ik} = \sum_j X_{ij}\otimes X_{jk}$ and $\varepsilon(X_{ik}) = \delta_{ik}$.

\noindent \textit{2.} $\newtilde{\Lambda}_q(n)$ is a right comodule algebra for $\newtilde{\mathrm{M}}_q(n)$ with $\rho(z_i) = \sum_j z_j\otimes X_{ji}$.
\end{myprop}

\begin{proof}[Proof (1)] We must check that the ideal $E$ is a coideal, i.e., that for all $I\in[n]^n$: (i) $\varepsilon(E_I)=0$ and (ii) $\Delta(E_I) \in E\otimes \K\langle X\rangle + \K\langle X\rangle \otimes E$. The verification of (i) is straightforward and is omitted. 

Working in $\K\langle X\rangle$, we have
\begin{eqnarray}
\nonumber \Delta(D_I) &=& \sum_{\pi\in\Sym_n} \Qlla{-}{\pi} \Delta \left(X_{\pi_1 i_1} \cdots X_{\pi_n i_n} \right)\\
\nonumber &=& \sum_{\pi\in\Sym_n} \Qlla{-}{\pi} \sum_{j_1,\ldots,j_n} X_{\pi_1 j_1} \cdots X_{\pi_n j_n} \otimes X_{j_1, i_1} \cdots X_{j_n i_n} \\
\nonumber &=& \sum_{j_1,\ldots,j_n} \left(\sum_{\pi\in\Sym_n} \Qlla{-}{\pi} X_{\pi_1 j_1} \cdots X_{\pi_n j_n} \right) \otimes X_{j_1, i_1} \cdots X_{j_n i_n} \\
\nonumber &=& \sum_{J=j_1,\ldots,j_n} D_J \otimes X_{j_1, i_1} \cdots X_{j_n i_n} \\
\nonumber \phantom{\Delta(D_I)}  &=& \sum_{J\not\in\Sym_n} E_J \otimes X_{j_1, i_1} \cdots X_{j_n i_n} + \sum_{J\in\Sym_n} E_J\otimes X_{j_1, i_1} \cdots X_{j_n i_n} \\
\nonumber &&\quad + \sum_{J\in\Sym_n} \Qlla{-}{J} D \otimes X_{j_1, i_1} \cdots X_{j_n i_n}  \\
\nonumber \phantom{\Delta(D_I)}   &=& \!\!\sum_{j_1,\ldots,j_n} E_J \otimes X_{j_1, i_1} \cdots X_{j_n i_n} \!+\!\! \sum_{J\in\Sym_n} D \otimes \Qlla{-}{J} X_{j_1, i_1} \!\cdots\! X_{j_n i_n} \\
\label{eq:grouplike} &=& \sum_{j_1,\ldots,j_n} E_J \otimes X_{j_1, i_1} \cdots X_{j_n i_n} + D \otimes D_I .
\end{eqnarray}
In particular, $D$ is grouplike in $\newtilde{\mathrm{M}}_q(n)$, and thus (ii) follows easily from (\ref{eq:grouplike}). 
\end{proof}

\begin{proof}[Proof (2)] It is clear that $\rho$ is a comodule map on $\newtilde{\Lambda}_q(n)_{(1)}$, the $\K$-subspace of $\newtilde{\Lambda}_q(n)$ spanned by the generators $z_i$. Extending $\rho$ by declaring it to be an algebra map, it is left to verify that it respects the relations $\mathcal E_I$ in $\newtilde{\Lambda}_q(n)$.

Suppose $I\in[n]^n$. We begin by computing $\rho(z_{i_1})\rho(z_{i_2})\cdots \rho(z_{i_n})$:
\begin{eqnarray}
\label{eq:step1} \rho(z_{i_1})\cdots \rho(z_{i_n}) &=& \sum_{j_1, \ldots, j_n} z_{j_1}\cdots z_{j_n} \otimes X_{j_1 i_1} \cdots X_{j_n i_n} \\
\label{eq:step2} &=& \sum_{J\in\Sym_n} \Qlla{-}{J} z_1\cdots z_n \otimes X_{j_1 i_1} \cdots X_{j_n i_n} \\
\label{eq:coinvariant} &=& z_1\cdots z_n \otimes D_I,
\end{eqnarray}
where, in moving from (\ref{eq:step1}) to (\ref{eq:step2}), we make use of the relations 
$\mathcal E_J$ to simplify the left tensor factors appearing in (\ref{eq:step1}). Now in order for $\rho$ to respect relation $\mathcal E_I$, we need $\rho(z_{i_1}\cdots z_{i_n} - \delta_{I,\Sym_n} \Qlla{-}{I} z_1 \cdots z_n)=0$. Using (\ref{eq:coinvariant}) above, it is clear that we need 
\[
z_1\cdots z_n \otimes (D_I - \delta_{I,\Sym_n}\Qlla{-}{I} D) =0,
\]
which holds by construction, namely, $E_I\equiv 0$ in $\newtilde{\mathrm{M}}_q(n)$.
\end{proof}

In the proof of Proposition \ref{thm:Mq-comodule}, we saw that $D$ is grouplike and $\rho(z_{1}z_{2}\cdots z_{n}) = z_1z_2\cdots z_n \otimes D$, cf. (\ref{eq:grouplike}) and (\ref{eq:coinvariant}) respectively. It follows immediately that $\newtilde{\mathrm{SL}}_q(n)$ is a bialgebra and that $\newtilde{\Lambda}_q(n)$ is also a $\newtilde{\mathrm{SL}}_q(n)$-comodule algebra. It remains to define a map $S$ on $\newtilde{\mathrm{SL}}_q(n)$ and show that it is a left but not a right antipode. We begin by building an amenable basis for $\newtilde{\mathrm{SL}}_q(n)$. 

Write the $n^n$ relations as reduction formulas, as illustrated earlier for $n=3$. Call the $n^n$ words appearing on the left side \emph{reducible} words of length $n$: 
\begin{equation}\label{eq:badwords} 
X_{n\,i_1} \cdots X_{2\,i_{n-1}} X_{1\,i_n}.
\end{equation}
A general word is \emph{irreducible} if it does not contain one of these as a subword. Now, for each $I=(i_1,\ldots,i_n)$, it is plain to see that the monomial in (\ref{eq:badwords}) occurs in one and only one relation defining $\newtilde{\mathrm{SL}}_q(n)$, namely the relation $D_I = \Qlla{-}{I} \delta_{I,\Sym_n}$. In particular, no reducible words appear on the right side of a reduction formula. Finally, there are no overlaps in the sense of the Diamond Lemma, i.e. no words $w=w_1w_2w_3$ with both $w_1w_2$ and $w_2w_3$ reducible. We conclude that the irreducible words (including the empty word ``$1$'') form a basis for $\newtilde{\mathrm{SL}}_q(n)$. 

Toward the goal of defining $S$, the key observation is that 
\begin{align}\tag{$\star$}
\begin{minipage}[c]{.85\textwidth}{\vspace{1.5ex}
\noindent \emph{if  $w = X_{i_1j_1}\cdots X_{i_r j_r}$ is an irreducible word, then 
$\Delta w = \sum_k w_k \otimes w'_k$  (before any reductions of the new words $w_k,w'_k$) 
 has the property that all $w_k$ are irreducible. 
\vspace{1.5ex}
}}\end{minipage}\end{align}
This follows from the form of the irreducible words. 

\begin{mynota}
Given any $m\times m$ matrix $A=(a_{ij})$, write $\qdet A$ for the expression
$\sum_{\pi\in\Sym_{m}} \Qlla{-}{\pi} a_{\pi_1 1}a_{\pi_2 2}\cdots a_{\pi_m m}$. Also, write $A^{ab}$ for the submatrix built from $A$ by deleting row $a$ and column $b$.
\end{mynota}

\begin{mydefi}\label{def:antipode} Define a set map $S:X\rightarrow \newtilde{\mathrm{SL}}_q(n)$ by the adjoint matrix, $S(X) = \mathrm{Adj}_q X$, i.e., $S(X_{ij}) = (-q)^{j-i} \qdet (X^{ji})$. Extend $S$ to an element of $\Hom(\newtilde{\mathrm{SL}}_q(n), \newtilde{\mathrm{SL}}_q(n))$ by following the scheme in \cite{NicTaf:1}, i.e., for irreducible words $w=X_{i_1 j_1}\cdots X_{i_r j_r}$ with $r>1$, put $S(w) = S(X_{i_r j_r})\cdots S(X_{i_1 j_1})$.
\end{mydefi}

\begin{mytheo}\label{thm:SLqn-properties} The map $S$ in Definition \ref{def:antipode} gives $\newtilde{\mathrm{SL}}_q(n)$ the structure of a left (but not a right) Hopf algebra.
\end{mytheo}

\begin{proof} Thanks to property ($\star$), it is a routine calculation to show that $m(S\otimes \II)\Delta(w) = \varepsilon(w)1$, once $S$ satisfies this property on the generators $X_{ij}$. 

Working in $\newtilde{\mathrm{M}}_q(n)$, we will show that $S(X)\cdot X = D\,\II_n$. This obviously proves half of the theorem. For the other half, we will show that $(\II\ast S)(X_{11}) = \sum_j X_{1j}S(X_{j1}) \neq 1$ in $\newtilde{\mathrm{SL}}_q(n)$. In both halves, we will make use of the following notation:

For any $j\in[n]$, denote the permutations of $[n]\setminus j$ by $\Sym_{[n]\setminus j}$. View the image $(\pi(1), \pi(2), \ldots, \widehat{\pi(j)},\ldots, \pi(n))$ of $\pi\in\Sym_{[n]\setminus j}$ as an ($n$\,$-1$)-tuple $(\pi_1, \ldots, \pi_{n-1})$. Extend $\ell$ to $\Sym_{[n]\setminus j}$ in the obvious manner.

\smallskip
\par\noindent\emph{Claim:} $\sum_j S(X_{ij})X_{ji'} = \delta_{ii'} D$.

\noindent Writing the left-hand side out carefully, we have 
\begin{eqnarray*}
\sum_j S(X_{ij})X_{ji'} %
&=& \sum_j (-q)^{j-i} \qdet (X^{ji}) X_{ji'} \\
&=& (-q)^{n-i}\sum_j  (-q)^{j-n} \Bigg(\sum_{\pi\in\Sym_{[n]\setminus j}} \Qlla{-}{\pi} \times \\
&& \qquad\qquad X_{\pi_1 1} \cdots X_{\pi_{i-1} i-1} X_{\pi_{i} i+1} \cdots X_{\pi_{n-1} n} \Bigg) X_{ji'} \\
&=& (-q)^{n-i} \sum_j \Bigg(\sum_{{\pi\in\Sym_{[n]\setminus j}}} \Qlla{-}{\pi_1,\ldots ,\pi_{n-1}, j} \times \\
&& \qquad\qquad X_{\pi_1 1} \cdots X_{\pi_{i-1} i-1}X_{\pi_{i} i+1}\cdots X_{\pi_{n-1} n}\Bigg) X_{ji'} \,, 
\end{eqnarray*}
since $\ell(\pi_1,\ldots,\pi_{n-1},j) = \ell(\pi_1,\ldots,\pi_{n-1}) + \ell(1,\ldots,\widehat{j},\ldots,n,j) = \ell(\pi) + n - j$. Also, since $\pi'=(\pi_1,\ldots,\pi_{n-1},j)$ runs over all permutations in $\Sym_n$, we conclude
\begin{eqnarray*}
\sum_j S(X_{ij})X_{ji'} &=& (-q)^{n-i} \sum_{{\pi'\in\Sym_n}} \Qlla{-}{\pi'} \times \\
&& \qquad X_{\pi'_1 1} \cdots X_{\pi'_{i-1} i-1} X_{\pi'_{i} i+1}\cdots X_{\pi'_{n-1} n} X_{\pi'_n i'} \\
&=& (-q)^{n-i} D_{(1,\ldots,\hat{i},\ldots, n,i')} \\
&=& (-q)^{n-i} \delta_{ii'} \Qlla{-}{1,\ldots,\widehat{i},\ldots,n,i}D = \delta_{ii'} D.
\end{eqnarray*}

\smallskip
\par\noindent\emph{Claim:} $\sum_j X_{1j}S(X_{j1}) \neq 1$. 

\noindent Here, the left-hand side is a linear combination of monomials of the form
\[
X_{1j}X_{\pi_1 1} \cdots X_{\pi_{j-1} j-1} X_{\pi_j j+1}\cdots X_{\pi_{n-1} n}.
\]
Each is a monomial of degree $n$ which is not of the form (\ref{eq:badwords}). In particular, they, together with $1$, form a linearly independent set in $\newtilde{\mathrm{SL}}_q(n)$.
\end{proof}

\def\cprime{$'$} \def\cprime{$'$} \def\cprime{$'$} \def\cprime{$'$}



\begin{thebibliography}{9}

\bibitem{Ber:1}
G.~M. Bergman, 
\newblock{`The diamond lemma for ring theory'},
\newblock {\em Adv. in Math.} {\bf 29} (1978) 2, 178--218.

\bibitem{FoaHan:2}
D. Foata, G.-N. Han,
\newblock{`A new proof of the {G}aroufalidis-{L}\^e-{Z}eilberger quantum {M}ac{M}ahon master theorem'},
\newblock{\em preprint} (2005).

\bibitem{GarLeZei:1}
S. Garoufalidis, T.~T. L\^e, D. Zeilberger,
\newblock{`The quantum {M}ac{M}ahon master theorem'}, 
\newblock{\em Proc. Natl. Acad. Sci. USA} (to appear).

\bibitem{GreNicTaf:1}
J.~A. Green, W.~D. Nichols, E.~J. Taft, 
\newblock{`Left {H}opf algebras'},
\newblock {\em J. Algebra} {\bf 65} (1980) 2, 399--411.

\bibitem{Jac:3}
N. Jacobson, 
\newblock{`Some remarks on one-sided inverses'},
\newblock {\em Proc. Amer. Math. Soc.} {\bf 1} (1950), 352--355.

\bibitem{NicTaf:1}
W.~D. Nichols, E.~J. Taft, 
\newblock{`The left antipodes of a left {H}opf algebra'},
\newblock in: {\em Algebraists' homage: papers in ring theory and related
  topics (New Haven, Conn., 1981)}, Vol.~13 of {\em Contemp. Math.}
\newblock Providence, R.I.: Amer. Math. Soc., 1982, pp. 363--368.

\bibitem{RodTaf:1}
S. Rodr{\'{\i}}guez-Romo, E.~J. Taft, 
\newblock{`Some quantum-like {H}opf algebras which remain noncommutative when {$q=1$}'},
\newblock {\em Lett. Math. Phys.} {\bf 61} (2002) 1, 41--50.

\bibitem{RodTaf:2}
S. Rodr{\'{\i}}guez-Romo, E.~J. Taft, 
\newblock{`A left quantum group'},
\newblock {\em J. Algebra} {\bf 286} (2005) 1, 154--160.

\bibitem{Tak:3}
M. Takeuchi, 
\newblock{`Matric bialgebras and quantum groups'},
\newblock {\em Israel J. Math.} {\bf 72} (1990) 1-2, 232--251.

\end{thebibliography}
\end{document}